\documentclass[final]{siamltex}
\usepackage{url,amsmath,graphicx,cite}
\usepackage[psamsfonts]{amssymb}
\usepackage[colorlinks, bookmarks=false]{hyperref}

\newcommand{\tr}{\operatorname{tr}}
\newcommand{\vecf}{\operatorname{vec}}

\title{Joint Singular Value Distribution of Two
Correlated Rectangular Gaussian Matrices and Its Application}

\author{Shuangquan Wang\thanks{Center for Wireless Communications and Signal
Processing Research (CWCSPR), Department of Electrical and Computer
Engineering, New Jersey Institute of Technology, Newark, NJ 07102
({\tt sw27@njit.edu}, {\tt ali.abdi@njit.edu}).}
        \and Ali Abdi${}^*$}

\begin{document}

\maketitle

\begin{abstract}
Let $\mathbf{H}=(h_{ij})$ and $\mathbf{G}=(g_{ij})$ be two $m\times
n$, $m\leq n$, random matrices, each with i.i.d complex zero-mean
unit-variance Gaussian entries, with correlation between any two
elements given by
$\mathbb{E}[h_{ij}g_{pq}^\star]=\rho\,\delta_{ip}\delta_{jq}$ such
that $|\rho|<1$, where ${}^\star$ denotes the complex conjugate and
$\delta_{ij}$ is the Kronecker delta. Assume $\{s_k\}_{k=1}^m$ and
$\{r_l\}_{l=1}^m$ are unordered singular values of $\mathbf{H}$ and
$\mathbf{G}$, respectively, and $s$ and $r$ are randomly selected
from $\{s_k\}_{k=1}^m$ and $\{r_l\}_{l=1}^m$, respectively. In this
paper, exact analytical closed-form expressions are derived for the
joint probability distribution function (PDF) of $\{s_k\}_{k=1}^m$
and $\{r_l\}_{l=1}^m$ using an Itzykson-Zuber-type integral, as well
as the joint marginal PDF of $s$ and $r$, by a bi-orthogonal
polynomial technique. These PDFs are of interest in multiple-input
multiple-output (MIMO) wireless communication channels and systems.
\end{abstract}



\pagestyle{myheadings} \thispagestyle{plain} \markboth{S. WANG AND
A. ABDI}{JOINT SINGULAR VALUE DISTRIBUTION OF CORRELATED MATRICES}

\section{Introduction} Random singular values have found numerous
applications such as hypothesis testing and principal component
analysis in statistics\cite{IEEE_sw27:MuirheadBook82}, nuclear
energy levels and level spacing in nuclear
physics\cite{IEEE_sw27:MehtaBook04}, and calculation of the
multiple-input multiple-output (MIMO) channel capacity in wireless
communications\cite{IEEE_sw27:Telatar99_MIMO_Cap}. The singular
value distribution of a \emph{single} Gaussian random matrix is
given in\cite{IEEE_sw27:Shen01_SVD_PDF}. However, the joint singular
value distribution of \emph{correlated} Gaussian random matrices
have received less attention so far, although it has important
applications in wireless MIMO communications, say, the second-order
statistics of the
\emph{eigen}-channels\cite{IEEE_sw27:Wang_EigenChannel_GlobeCom06}
and instantaneous mutual
information\cite{IEEE_sw27:Wang05_LCR_AFD_SISO_Capacity,
IEEE_sw27:Wang05_IT,IEEE_sw27:NanZhang05}.

To the best of our knowledge, correlated random matrices have been
studied to some extent\cite{IEEE_sw27:MehtaBook04,
IEEE_sw27:Mehta94_Two_Coupled_Hermitian_Matrices,
IEEE_sw27:Mehta98_Coupled_Hermitian_Matrix_Chain}, where only
Hermitian matrices were considered. Different from
\cite{IEEE_sw27:MehtaBook04,
IEEE_sw27:Mehta94_Two_Coupled_Hermitian_Matrices,
IEEE_sw27:Mehta98_Coupled_Hermitian_Matrix_Chain}, we consider the
situation where the elements, with the same indices, of the two
rectangular complex Gaussian random matrices are correlated by a
\emph{complex} number, and derive exact analytical closed-form
expressions for the joint PDF of their singular values.

This paper is organized as follows. Section \ref{sec:problem_desc}
introduces the two rectangular complex Gaussian random matrices. The
joint PDF of singular values are studied in section
\ref{sec:JPDF_Singular_Values} using an Itzykson-Zuber-type
integral. The joint marginal PDF of singular values is derived in
section \ref{sec:Marginal_JPDF} and its application to wireless MIMO
communications is presented in section \ref{sec:Applications}.
Finally, concluding remarks are summarized in section
\ref{sec:Conclusion}.

{\it Notation}: $\cdot^{\dag}$ is reserved for matrix Hermitian,
$\cdot^T$ for matrix transpose, $\cdot^{\star}$ for complex
conjugate, $\tr[\cdot]$ for the trace of a matrix, $\jmath$ for
$\sqrt{-1}$, $\mathbb{E}[\cdot]$ for mathematical expectation,
$\mathbf{I}_m$ for the $m\times m$ identity matrix, $\otimes$ for
the Kronecker product, and $\Re[\cdot]$ and $\Im[\cdot]$ for the
real and imaginary parts of a complex number, respectively. In
addition, $\diag(\mathbf{s})$ denotes a diagonal matrix with
$\mathbf{s}$ on the main diagonal, $t\!\!\in\!\![m,n]$ implies that
$t$, $m$, and $n$ are all integers such that $m\leq t\leq n$ with
$m\leq n$, and $\det\left|x_{kl}\right|$ is the determinant of the
matrix, where $x_{kl}$ resides on the $k^\mathrm{th}$ row and
$l^\mathrm{th}$ column. Moreover, lower-case bold letters represent
row vectors, whereas upper-case bold letters are used for matrices.
Finally $\mathcal{C\,\!N}$ means complex normal, and $\vecf(\cdot)$
stacks all the columns of its matrix argument into one tall column
vector.

\section{Problem Description}\label{sec:problem_desc}

There are two $m\times n$ random matrices $\mathbf{H}=(h_{ij})$ and
$\mathbf{G}=(g_{ij})$, $i\in[1,m]$, $j\in[1,n]$, each with i.i.d
complex zero-mean unit-variance Gaussian entries, i.e.,
$\mathbb{E}[h_{ij}]=\mathbb{E}[g_{ij}]=0, \forall i, j$,
$\mathbb{E}[h_{ij}h_{pq}^\star]=\mathbb{E}[g_{ij}g_{pq}^\star]
=\delta_{ip}\delta_{jq}$, where the Kronecker symbol $\delta_{ij}$
is $1$ or $0$ when $i=j$ or $i\neq j$. Therefore $\mathbf{H},
\mathbf{G}\thicksim \mathcal{C\,\!N}(\mathbf{0}, \mathbf{I}_{mn})$.
Moreover, the correlation among the two random matrices is given by
\begin{equation}\label{eq:autocorrelation}
\mathbb{E}[h_{ij}g_{pq}^\star]=\rho\,\delta_{ip}\delta_{jq}, \quad
\forall i,j,p,q,
\end{equation} where $\rho=|\rho|e^{\jmath\theta}$ is a complex
number with $|\rho|<1$.

Without loss of generality, we assume $m\leq n$ and set $\nu=n-m$.
Based on the singular value decomposition (SVD), $\mathbf{H}$ and
$\mathbf{G}$ can be, respectively, diagonalized
as\cite{IEEE_sw27:HuaBook63}
\begin{align}
\mathbf{H}&=\mathbf{U}\mathbf{S}\mathbf{V}^\dag,\label{eq:diagonalization_H}\\
\mathbf{G}&=\widetilde{\mathbf{U}}\mathbf{R}
\widetilde{\mathbf{V}}^\dag,\label{eq:diagonalization_G}
\end{align} where
$\mathbf{S}=\begin{bmatrix}\diag(\mathbf{s})\
\mathbf{0}\end{bmatrix}$ and
$\mathbf{R}=\begin{bmatrix}\diag(\mathbf{r}) \
\mathbf{0}\end{bmatrix}$ with $\mathbf{s}=[s_1, s_2, \cdots, s_m]$
and $\mathbf{r}=[r_1, r_2, \cdots, r_m]$, respectively.

We assume that the singular values of $\mathbf{G}$, $r_1, r_2,
\cdots, r_m$, are unordered and the singular values of $\mathbf{H}$,
$s_1, s_2, \cdots, s_m$, are also unordered. Now we would like to
know the joint PDF of $\{r_l\}_{l=1}^m$ and $\{s_l\}_{l=1}^m$.
Moreover, with $r$ randomly selected from $r_1, r_2, \cdots, r_m$,
and $s$ randomly selected from $s_1, s_2, \cdots, s_m$, it is of
interest to derive the joint PDF of $r$ and $s$ as well. These two
PDFs are derived in Section \ref{sec:JPDF_Singular_Values} and
\ref{sec:Marginal_JPDF}, respectively.

\section{Joint PDF of $\{s_l\}_{l=1}^m$ and
$\{r_l\}_{l=1}^m$}\label{sec:JPDF_Singular_Values}

\begin{lemma}[{\rm Joint PDF of $\mathbf{H}$ and
$\mathbf{G}$}] \label{lem:jpdf_H_G} For two correlated rectangular
complex Gaussian random matrices, $\mathbf{H}, \mathbf{G}\thicksim
\mathcal{C\,\!N}(\mathbf{0}, \mathbf{I}_{mn})$, with the correlation
between $\mathbf{H}$ and $\mathbf{G}$ given by
(\ref{eq:autocorrelation}), the joint PDF of $\mathbf{H}$ and
$\mathbf{G}$ is given by
\begin{equation}\label{eq:jpdf_H_G}
p(\mathbf{H},
\mathbf{G})=\frac{1}{\pi^{2mn}\left(1-|\rho|^2\right)^{mn}}
\exp\left[-\frac{\tr\!\left(\mathbf{H}\mathbf{H}^\dag+
\mathbf{G}\mathbf{G}^\dag-\rho^\star\mathbf{H}\mathbf{G}^\dag-
\rho\mathbf{G}\mathbf{H}^\dag\right)}{1-|\rho|^2}\right].
\end{equation}
\end{lemma}
\begin{proof}
We set $\mathbf{h}=\vecf(\mathbf{H})$,
$\mathbf{g}=\vecf(\mathbf{G})$, and $\mathbf{x}=[\mathbf{h}^T\
\mathbf{g}^T]^T$. Based on $\mathbf{H}, \mathbf{G}\thicksim
\mathcal{C\,\!N}(\mathbf{0}, \mathbf{I}_{mn})$ and
(\ref{eq:autocorrelation}), we have the mean and covariance matrix
of $\mathbf{x}$ as $\mathbb{E}[\mathbf{x}]=\mathbf{0}$ and
$\Sigma_\mathbf{x}=\Sigma_\tau\otimes\mathbf{I}_{mn}$ with
$\Sigma_\tau=\left[\begin{smallmatrix}1& \rho\\\rho^\star&
1\end{smallmatrix}\right]$, respectively. Therefore the PDF of
$\mathbf{x}$ is given by\cite{IEEE_sw27:James64_MatrixVariate}
\begin{equation}\label{eq:pdf_x}
p(\mathbf{x})=\frac{1}{\pi^{2mn}\det|\Sigma_\mathbf{x}|}
\exp\left(-\mathbf{x}^\dag\Sigma_\mathbf{x}^{-1}\mathbf{x}\right),
\end{equation} where $\det|\Sigma_\mathbf{x}|=
\left(\det|\Sigma_\tau|\right)^{mn}=\left(1-|\rho|^2\right)^{mn}$.

With
$\Sigma_\tau^{-1}=\frac{1}{1-|\rho|^2}\left[\begin{smallmatrix}1&
-\rho\\-\rho^\star& 1\end{smallmatrix}\right]$, we obtain
$\Sigma_\mathbf{x}^{-1}=\Sigma_\tau^{-1}\otimes\mathbf{I}_{mn}=
\frac{1}{1-|\rho|^2}\left[\begin{smallmatrix}\mathbf{I}_{mn}&
-\rho\mathbf{I}_{mn}\\-\rho^\star\mathbf{I}_{mn}&
\mathbf{I}_{mn}\end{smallmatrix}\right]$. Therefore
$\mathbf{x}^\dag\Sigma_\mathbf{x}^{-1}\mathbf{x}$ in
(\ref{eq:pdf_x}) can be rewritten as
\begin{equation}\label{eq:trace}
\begin{split}
\mathbf{x}^\dag\Sigma_\mathbf{x}^{-1}\mathbf{x}
&=\tr\left(\Sigma_\mathbf{x}^{-1}\mathbf{x}\mathbf{x}^\dag\right)
=\tr\left(\frac{1}{1-|\rho|^2}\left[\begin{smallmatrix}\mathbf{I}_{mn}&
-\rho\mathbf{I}_{mn}\\-\rho^\star\mathbf{I}_{mn}&
\mathbf{I}_{mn}\end{smallmatrix}\right]\left[\begin{smallmatrix}
\mathbf{h}\mathbf{h}^\dag& \mathbf{h}\mathbf{g}^\dag\\
\mathbf{g}\mathbf{h}^\dag&
\mathbf{g}\mathbf{g}^\dag\end{smallmatrix}\right]\right),\\
&=\frac{\tr\left(\mathbf{h}\mathbf{h}^\dag
+\mathbf{g}\mathbf{g}^\dag-\rho^\star\mathbf{h}\mathbf{g}^\dag
-\rho\mathbf{g}\mathbf{h}^\dag\right)}{1-|\rho|^2}
=\frac{\tr\!\left(\mathbf{H}\mathbf{H}^\dag+
\mathbf{G}\mathbf{G}^\dag-\rho^\star\mathbf{H}\mathbf{G}^\dag-
\rho\mathbf{G}\mathbf{H}^\dag\right)}{1-|\rho|^2},
\end{split}
\end{equation} where $\tr\left(\mathbf{A}\mathbf{B}^\dag\right)=
\vecf(\mathbf{B})^\dag\vecf(\mathbf{A}) =\tr\left[\vecf(\mathbf{A})
\vecf(\mathbf{B})^\dag\right]$\cite{IEEE_sw27:GuptaBook99} is used
in the last ``='' of (\ref{eq:trace}). Substitution of
(\ref{eq:trace}) into (\ref{eq:pdf_x}) leads to (\ref{eq:jpdf_H_G}).
\end{proof}

From (\ref{eq:diagonalization_H}), we know that the unitary matrix
pair $(\mathbf{U}, \mathbf{V})$ parameterizes the coset space
$\mathcal{U}(m)\times\mathcal{U}(n)/\left[\mathcal{U}(1)\right]^m$,
where $\mathcal{U}(p)$ is the unitary group of order $p$, and the
integration measure,
$d[\mathbf{H}]=\prod_{i=1}^m\prod_{j=1}^nd\left[\Re{h_{ij}}\right]
d\left[\Im{h_{ij}}\right]$, can be represented
by\cite{IEEE_sw27:Jackson96}
\begin{equation}\label{eq:dH}
d[\mathbf{H}]=\Omega J(\mathbf{s})d[\mathbf{s}]d\mu(\mathbf{U},
\mathbf{V}),
\end{equation} where $J(\mathbf{s})=\triangle^2(\mathbf{s}^2)
\prod_{k=1}^ms_k^{2\nu+1}$ with the $m$-dimensional Vandermonde
determinant $\triangle(\mathbf{s}^2)=\det\left|s_k^{2(l-1)}\right|
=\prod_{k>l}(s^2_k-s^2_l)$ and
$\triangle^2(\cdot)=\left[\triangle(\cdot)\right]^2$,
$d[\mathbf{s}]=\prod_{l=1}^m ds_l$, $d\mu(\mathbf{U}, \mathbf{V})$
is the Haar measure of $\mathcal{U}(m)\times\mathcal{U}(n)/
\left[\mathcal{U}(1)\right]^m$\cite{IEEE_sw27:Jackson96}, and the
constant $\Omega$ is given
by\cite{IEEE_sw27:Nagao91,IEEE_sw27:Jackson96}
\begin{equation}\label{eq:Omega}
\Omega=\frac{2^m\pi^{mn}}{\prod_{j=1}^mj!(j+\nu-1)!}
=\frac{2^m\pi^{mn}}{m!\prod_{j=0}^{m-1}j!(j+\nu)!}.
\end{equation}

Similarly, we have
\begin{equation}\label{eq:dG}
d[\mathbf{G}]=\Omega
J(\mathbf{r})d[\mathbf{r}]d\mu(\widetilde{\mathbf{U}},
\widetilde{\mathbf{V}}),
\end{equation} where $J(\mathbf{r})=\triangle^2(\mathbf{r}^2)
\prod_{k=1}^mr_k^{2\nu+1}$ with the $m$-dimensional Vandermonde
determinant $\triangle(\mathbf{r}^2)=\det \left|r_k^{2(l-1)}\right|
=\prod_{k>l}(r^2_k-r^2_l)$ and $d[\mathbf{r}]=\prod_{l=1}^m dr_l$.

In order to obtain the joint PDF of $\{r_l\}_{l=1}^m$ and
$\{s_l\}_{l=1}^m$, we need the following proposition.
\begin{proposition}[{\rm An Itzykson-Zuber-type
integral\cite[(31)]{IEEE_sw27:Jackson96}}]
\label{pro:Itzykson-Zuber-Integral}
\begin{equation}
\begin{split}
&\hspace{1em}\int d\mu(\mathbf{U},\mathbf{V})\exp\left\{-\frac{\tr
\left[(\mathbf{H}-\mathbf{G})
(\mathbf{H}-\mathbf{G})^\dag\right]}{t}\right\}\\
&=\frac{2^m\pi^{mn}t^{mn-m}
\det\left|\exp\left(-\frac{s^2_k+r^2_l}{t}\right)
I_\nu\!\!\left(\frac{2s_kr_l}{t}\right)\right|}
{m!\Omega\triangle(\mathbf{s}^2)
\triangle(\mathbf{r}^2)\prod_{k=1}^m (s_kr_k)^\nu},
\end{split}
\end{equation} where $\Omega$ is given by (\ref{eq:Omega}) and
$I_k(z)=\frac{1}{\pi}\!\int_0^\pi e^{z\cos\theta}
\cos(k\theta)\text{d}\theta$ is the $k^\text{th}$ order modified
Bessel function of the first kind.
\end{proposition}

\begin{theorem}\label{th:jpdf_singular_values}
The joint PDF of the singular values of $\mathbf{H}$ and
$\mathbf{G}$ is given by
\begin{equation}\label{eq:jpdf_singular_values}
p(\mathbf{s}, \mathbf{r})
=\frac{\exp\left(-\frac{\sum_{k=1}^ms^2_k+r^2_k}
{1-|\rho|^2}\right)\triangle(\mathbf{s}^2)
\triangle(\mathbf{r}^2)\prod_{k=1}^m (s_kr_k)^{\nu+1} \det\left|
I_\nu\!\!\left(\frac{2|\rho|s_kr_l}{1-|\rho|^2}\right)\right|}
{2^{-2m}m!m!\prod_{j=0}^{m-1}j!(j+\nu)!|\rho|^{mn-m}(1-|\rho|^2)^m}.
\end{equation}
\end{theorem}
\begin{proof}
By combining (\ref{eq:jpdf_H_G}) with (\ref{eq:dH}) and
(\ref{eq:dG}), we obtain
\begin{equation}\label{eq:jpdf_singular_values_derivation}
p(\mathbf{s}, \mathbf{r})= \frac{\Omega^2
J(\mathbf{s})J(\mathbf{r})}{\pi^{2mn}(1-|\rho|^2)^{mn}}
\Phi(\mathbf{s}, \mathbf{r}),
\end{equation} where
\begin{equation}\label{eq:Phi_s_r}
\begin{split}
\Phi(\mathbf{s}, \mathbf{r})\!&=\!\int\!\!
d\mu(\widetilde{\mathbf{U}},\widetilde{\mathbf{V}})\!\int\!\!
d\mu(\mathbf{U},\mathbf{V})\exp\left[-\frac{\tr\!\left(
\mathbf{H}\mathbf{H}^\dag
\!+\!\mathbf{G}\mathbf{G}^\dag\!-\!\rho^\star\mathbf{H}\mathbf{G}^\dag\!-\!
\rho\mathbf{G}\mathbf{H}^\dag\right)}{1-|\rho|^2}\right],\\
\!&=\!\int\!\!
d\mu(\widetilde{\mathbf{U}},\widetilde{\mathbf{V}})\!\!\int\!\!
d\mu(\mathbf{U},\mathbf{V})\exp\left\{\!-\frac{\tr
\left[(\mathbf{H}\!-\!\rho\mathbf{G})
(\mathbf{H}\!-\!\rho\mathbf{G})^\dag\right]}{1-|\rho|^2}\!-\!
\tr(\mathbf{G}\mathbf{G}^\dag)\!\right\}\!,\\
\!&=\!\int\!\! d\mu(\widetilde{\mathbf{U}},\widetilde{\mathbf{V}})
\,e^{-\tr(\mathbf{G}\mathbf{G}^\dag)}\!\int\!\!
d\mu(\mathbf{U},\mathbf{V})\exp\left\{\!-\frac{\tr
\left[(\mathbf{H}-\rho\mathbf{G})
(\mathbf{H}-\rho\mathbf{G})^\dag\right]}{1-|\rho|^2}\!\right\},\\
\!&=\!\int\!\! d\mu(\widetilde{\mathbf{U}},\widetilde{\mathbf{V}})
\frac{e^{-\sum_{k=1}^mr_k^2}(1-|\rho|^2)^{mn-m}
\det\left|e^{-\frac{s^2_k+|\rho|^2r^2_l}{1-|\rho|^2}}
I_\nu\!\!\left(\frac{2|\rho|s_kr_l}{1-|\rho|^2}\right)\right|}
{2^{-m}m!\pi^{-mn}\Omega\triangle(\mathbf{s}^2)
\triangle(|\rho|^2\mathbf{r}^2)\prod_{k=1}^m
(|\rho|s_kr_k)^\nu},\\
\!&=\!\frac{(1-|\rho|^2)^{mn-m}\exp\left(-\frac{\sum_{k=1}^ms^2_k+r^2_k}
{1-|\rho|^2}\right) \det\left|
I_\nu\!\!\left(\frac{2|\rho|s_kr_l}{1-|\rho|^2}\right)\right|}
{2^{-m}m!\pi^{-mn}\Omega|\rho|^{mn-m}\triangle(\mathbf{s}^2)
\triangle(\mathbf{r}^2)\prod_{k=1}^m (s_kr_k)^\nu}.
\end{split}
\end{equation} Derivation of the second and third lines of
(\ref{eq:Phi_s_r}) are straightforward. The fourth line comes from %(\ref{eq:drhoG}) and
\begin{equation}
\rho\mathbf{G}=\widehat{\mathbf{U}}\widehat{\mathbf{R}}
\widehat{\mathbf{V}}^\dag
\end{equation} with
$\widehat{\mathbf{R}}=|\rho|\mathbf{R}$, and Proposition
\ref{pro:Itzykson-Zuber-Integral} with the replacements
$t\rightarrow1-|\rho|^2$ and $\mathbf{G}\rightarrow\rho\mathbf{G}$.
The last line is based on the convention that $\int\!
d\mu(\widetilde{\mathbf{U}},\widetilde{\mathbf{V}})=1$\cite{IEEE_sw27:Jackson96}.
Plugging (\ref{eq:Omega}) and the last line of (\ref{eq:Phi_s_r})
into (\ref{eq:jpdf_singular_values_derivation}), we obtain
(\ref{eq:jpdf_singular_values}).
\end{proof}

By relating the eigenvalues of $\mathbf{G}\mathbf{G}^\dag$ to the
singular values of $\mathbf{G}$ through $\alpha_l=r_l^2$,
$l\in[1,m]$, and the eigenvalues of $\mathbf{H}\mathbf{H}^\dag$ to
the singular values of $\mathbf{H}$ through $\beta_l=s_l^2$,
$l\in[1,m]$, we can derive the joint PDF of
$\mbox{\boldmath{$\alpha$}}=\left[\alpha_1, \alpha_2, \cdots,
\alpha_m\right]$ and $\mbox{\boldmath{$\beta$}}=\left[\beta_1,
\beta_2, \cdots, \beta_m\right]$, presented in the following
corollary.
\begin{corollary}\label{coro:jpdf_eigenvalues}
The joint PDF of the unordered eigenvalues of
$\mathbf{H}\mathbf{H}^\dag$ and $\mathbf{G}\mathbf{G}^\dag$ is
\begin{equation}\label{eq:jpdf_eig_values}
p(\mbox{\boldmath{$\beta$}}, \mbox{\boldmath{$\alpha$}})=\frac{\exp
\left(-\frac{\sum_{k=1}^m\beta_k+\alpha_k}
{1-|\rho|^2}\right)\triangle(\mbox{\boldmath{$\beta$}})
\triangle(\mbox{\boldmath{$\alpha$}})\prod_{k=1}^m
(\sqrt{\beta_k\alpha_k})^\nu \det\left|
I_\nu\!\!\left(\frac{2|\rho|\sqrt{\beta_k\alpha_l}}{1-|\rho|^2}\right)\right|}
{m!m!\prod_{j=0}^{m-1}j!(j+\nu)!|\rho|^{mn-m}(1-|\rho|^2)^m},
\end{equation} where $m$-dimensional Vandermonde determinants are defined by
$\triangle(\mbox{\boldmath{$\beta$}})=\det
\left|\beta_k^{l-1}\right|=\prod_{k>l}(\beta_k-\beta_l)$ and
$\triangle(\mbox{\boldmath{$\alpha$}})=\det
\left|\alpha_k^{l-1}\right|=\prod_{k>l}(\alpha_k-\alpha_l)$.
\end{corollary}
\begin{proof}
It is straightforward to obtain (\ref{eq:jpdf_eig_values}) from
(\ref{eq:jpdf_singular_values}) by $2m$ one-to-one nonlinear
mappings.
\end{proof}
\section{Joint Marginal PDF}\label{sec:Marginal_JPDF}
In this section, with $\beta=s^2$ and $\alpha=r^2$, we calculate the
joint marginal PDF of $\beta$ and $\alpha$, $p(\beta, \alpha)$,
using the techniques and results presented in
\cite{IEEE_sw27:Mehta94_Two_Coupled_Hermitian_Matrices,
IEEE_sw27:Mehta98_Coupled_Hermitian_Matrix_Chain}. Then the joint
PDF of $s$ and $r$, $p(s,r)$, is easily derived.

If the polynomials $P_k(\beta)$ and $Q_l(\alpha)$, satisfy $\int
w(\beta, \alpha)P_k(\beta)Q_l(\alpha)d\beta d\alpha=\delta_{kl}$,
then we call $P_k(\beta)$ and $Q_l(\alpha)$ as bi-orthogonal
polynomials, associated with the weight function $w(\beta,
\alpha)$\cite{IEEE_sw27:MehtaBook04}. With this definition, we have
the following Lemma.
\begin{lemma}\label{lem:jpdf_by_bi-orthonomal_weight_function}
There exist bi-polynomials, $P_k(\beta)$ and $Q_l(\alpha)$, and a
weight function, $w(\beta, \alpha)$, which reduce
(\ref{eq:jpdf_eig_values}) to the following form
\begin{equation}\label{eq:jpdf_by_bi-orthonomal_weight_function}
p(\mbox{\boldmath{$\beta$}},
\mbox{\boldmath{$\alpha$}})=C_1\,\det|P_{k-1}(\beta_l)|
\det|w(\beta_k, \alpha_l)| \det|Q_{k-1}(\alpha_l)|,
\end{equation} where $C_1$ is a normalization constant. % where $C=m!m!$
\end{lemma}
\begin{proof}
In this paper, $\nu$ is a non-negative integer. Using the
Hille-Hardy formula\cite[pp.~185, (46)]{IEEE_sw27:BeckmannBook73}
\begin{equation}\label{eq:Hille-Hardy-Formula}
\sum_{k=0}^\infty\frac{k!z^k}{(k+\nu)!}L_k^\nu(x)L_k^\nu(y)
=\frac{(xyz)^{-\frac{\nu}{2}}}{1-z}\exp\left(-z\frac{x+y}{1-z}\right)
I_\nu\!\!\left(\frac{2\sqrt{xyz}}{1-z}\right), |z|<1,
\end{equation} with $L_k^\nu(x)=\frac{1}{k!}
e^xx^{-\nu}\frac{d^k}{dx^k}(e^{-x}x^{k+\nu})$ as the associated
Laguerre polynomial, we can rewrite (\ref{eq:jpdf_eig_values}) as
\begin{equation}\label{eq:jpdf_eig_values_Laguerre}
p(\mbox{\boldmath{$\beta$}},
\mbox{\boldmath{$\alpha$}})=\frac{\triangle(\mbox{\boldmath{$\beta$}})
\triangle(\mbox{\boldmath{$\alpha$}})\det\left| \beta_k^\nu
e^{-\beta_k}\alpha_l^\nu e^{-\alpha_l}\sum_{j=0}^\infty
\frac{j!|\rho|^{2j}L_j^\nu(\beta_k)L_j^\nu(\alpha_l)}{(j+\nu)!}\right|}
{m!m!\prod_{j=0}^{m-1}j!(j+\nu)!|\rho|^{m(m-1)}}.
\end{equation}

We set the weight function, $w(\beta, \alpha)$, as
\begin{equation}\label{eq:weight_function_w_beta_alpha}
\begin{split}
w(\beta, \alpha)&=\beta^\nu\alpha^\nu
e^{-(\beta+\alpha)}\sum_{j=0}^\infty
\frac{j!|\rho|^{2j}L_j^\nu(\beta)L_j^\nu(\alpha)}{(j+\nu)!},\\
&=\frac{(\beta\alpha)^\frac{\nu}{2}
e^{-\frac{\beta+\alpha}{1-|\rho|^2}}
I_\nu\!\!\left(\frac{2|\rho|\sqrt{\beta\alpha}}
{1-|\rho|^2}\right)}{(1-|\rho|^2)|\rho|^\nu}.
\end{split}
\end{equation} It is easy to check that
the corresponding bi-orthogonal polynomials are given by
\begin{align}
  P_k(\beta) &=\sqrt{\frac{k!}{(k+\nu)!}}|\rho|^{-k}L_k^\nu(\beta),
  \label{eq:Bi-orthogonal_Poly_P}\\
  Q_l(\alpha)
  &=\sqrt{\frac{l!}{(l+\nu)!}}|\rho|^{-l}L_l^\nu(\alpha),
  \label{eq:Bi-orthogonal_Poly_Q}
\end{align} using the following integral
equality\cite[pp.~267, 7.414.3]{IEEE_sw27:BeckmannBook73}
\begin{equation}\label{eq:Orthorgonality_Laguerre_Poly}
\int_0^\infty e^{-x}x^\nu L_k^\nu(x)L_l^\nu(x)
=\frac{(k+\nu)!}{k!}\delta_{kl}
\end{equation}

Moreover, by the addition of multiples of rows of lower order which
do not change the determinant of the Vandermonde matrix, then each
of the rows can be expressed in terms of orthogonal polynomials with
respect to the weight function $w(\beta, \alpha)$. Therefore two
$m$-dimensional Vandermonde determinants,
$\triangle(\mbox{\boldmath{$\beta$}})$ and
$\triangle(\mbox{\boldmath{$\alpha$}})$, can be represented as
\begin{align}
\triangle(\mbox{\boldmath{$\beta$}})
&=\det\left|\beta_k^{l-1}\right|
=C_2\det\left|P_{k-1}(\beta_l)\right|,\label{eq:Vandermode_det_P}\\
\triangle(\mbox{\boldmath{$\alpha$}})
&=\det\left|\alpha_k^{l-1}\right|
=C_3\det\left|Q_{k-1}(\alpha_l)\right|,\label{eq:Vandermode_det_Q}
\end{align} where we use the fact that the matrix transpose does not
change the determinant, i.e., $\det\left|P_{l-1}(\beta_k)\right|
=\det\left|P_{k-1}(\beta_l)\right|$ and
$\det\left|Q_{l-1}(\beta_k)\right|
=\det\left|Q_{k-1}(\beta_l)\right|$.

The coefficient of $x^k$ in $L_k^\nu(x)$ is $\frac{(-1)^k}{k!}$,
then the coefficient of $x^k$ in $P_k(x)$ is
$(-1)^k|\rho|^{-k}\frac{1}{\sqrt{k!(k+\nu)!}}$, therefore we have
$C_2=\prod_{j=0}^{m-1}\!(-1)^j|\rho|^j\sqrt{j!(j\!+\!\nu)!}=
(-1)^{\frac{m(m-1)}{2}}\\\times\sqrt{|\rho|^{m(m-1)}\!
\prod_{j=0}^{m-1}\!j!(j\!+\!\nu)!}$, obtained by plugging
(\ref{eq:Bi-orthogonal_Poly_P}) into (\ref{eq:Vandermode_det_P}).
Similarly, substitution of (\ref{eq:Bi-orthogonal_Poly_Q}) into
(\ref{eq:Vandermode_det_Q}) gives $C_3=C_2$. Now the product of
(\ref{eq:Vandermode_det_P}) and (\ref{eq:Vandermode_det_Q}) results
in
\begin{equation}\label{eq:Vandermode_det_relationship}
\triangle(\mbox{\boldmath{$\beta$}})\triangle(\mbox{\boldmath{$\alpha$}})
=|\rho|^{m(m-1)}\prod_{j=0}^{m-1}j!(j+\nu)!\det|P_{k-1}(\beta_l)|
\det|Q_{k-1}(\alpha_l)|.
\end{equation}

Based on (\ref{eq:weight_function_w_beta_alpha}) and
(\ref{eq:Vandermode_det_relationship}), one can see that
(\ref{eq:jpdf_eig_values_Laguerre}) is equal to
(\ref{eq:jpdf_by_bi-orthonomal_weight_function}) with
$C_1=\frac{1}{m!m!}$.
\end{proof}

\begin{theorem}\label{th:jpdf_beta_alpha}
The joint PDF of $\beta$ and $\alpha$ is given by
\begin{multline}\label{eq:jpdf_beta_alpha}
p(\beta, \alpha)=\frac{(\beta\alpha)^\frac{\nu}{2}
e^{-\frac{\beta+\alpha}{1-|\rho|^2}}I_\nu\!\!\left(\frac{2|\rho|\sqrt{\beta\alpha}}
{1-|\rho|^2}\right)}{m^2(1-|\rho|^2)|\rho|^\nu}
\sum_{k=0}^{m-1}\frac{k!}{(k+\nu)!}\frac{L_k^\nu(\beta)L_k^\nu(\alpha)}
{|\rho|^{2k}}\\
+\frac{(\beta\alpha)^\nu e^{-(\beta+\alpha)}}{m^2}\sum_{0\leq
k<l}^{m-1}\bigg\{\frac{k!l!}{(k+\nu)!(l+\nu)!}
\Big\{\left[L_k^\nu(\beta)L_l^\nu(\alpha)\right]^2
+\left[L_l^\nu(\beta)L_k^\nu(\alpha)\right]^2\\
-\left[|\rho|^{2(l-k)}+|\rho|^{2(k-l)}\right]
L_k^\nu(\beta)L_l^\nu(\beta)L_k^\nu(\alpha)L_l^\nu(\alpha)\Big\}\bigg\}.
\end{multline}
\end{theorem}
\begin{proof}
Based on Lemma \ref{lem:jpdf_by_bi-orthonomal_weight_function}, and
the results presented in
\cite[(3.7)]{IEEE_sw27:Mehta94_Two_Coupled_Hermitian_Matrices}
\cite{IEEE_sw27:Mehta98_Coupled_Hermitian_Matrix_Chain}, $p(\beta,
\alpha)$ can be expressed as
\begin{equation}\label{eq:jpdf_beta_alpha_general}
m^2\,p(\beta, \alpha)=w(\beta,
\alpha)\sum_{k=0}^{m-1}P_k(\beta)Q_k(\alpha)+\sum_{0\leq
k<l}^{m-1}\det\!\left|\!\!\!
                \begin{array}{cc}
                  P_k(\beta) & \!\!\overline{P}_k(\alpha) \\
                  P_l(\beta) & \!\!\overline{P}_l(\alpha) \\
                \end{array}
              \!\!\!\right|
              \!\det\!\left|\!\!\!
                \begin{array}{cc}
                  \overline{Q}_k(\beta) & \!\!Q_k(\alpha) \\
                  \overline{Q}_l(\beta) & \!\!Q_l(\alpha) \\
                \end{array}
              \!\!\!\right|,
\end{equation} where $P_k(x)$ and $Q_k(x)$ are defined in (\ref{eq:Bi-orthogonal_Poly_P})
and (\ref{eq:Bi-orthogonal_Poly_Q}), respectively, the weight
function is presented in (\ref{eq:weight_function_w_beta_alpha}),
and $\overline{P}_k(\alpha)$ and $\overline{Q}_l(\beta)$ are
similarly defined as
\cite{IEEE_sw27:Mehta94_Two_Coupled_Hermitian_Matrices}
\begin{align}%\label{eq:Weighted_Bi-orthogonal_Poly}
  \overline{P}_k(\alpha) &=\int P_k(\beta)w(\beta, \alpha)d\beta
  =\sqrt{\frac{k!}{(k+\nu)!}}\alpha^\nu e^{-\alpha}|\rho|^kL_k^\nu(\alpha),
  \label{eq:Weighted_Bi-orthogonal_Poly_P}\\
  \overline{Q}_l(\beta) &=\int Q_l(\alpha)w(\beta, \alpha)d\alpha
  =\sqrt{\frac{l!}{(l+\nu)!}}\beta^\nu
  e^{-\beta}|\rho|^lL_l^\nu(\beta).\label{eq:Weighted_Bi-orthogonal_Poly_Q}
\end{align}

Plugging (\ref{eq:weight_function_w_beta_alpha}),
(\ref{eq:Bi-orthogonal_Poly_P}), (\ref{eq:Bi-orthogonal_Poly_Q}),
(\ref{eq:Weighted_Bi-orthogonal_Poly_P}) and
(\ref{eq:Weighted_Bi-orthogonal_Poly_Q}) into
(\ref{eq:jpdf_beta_alpha_general}), we arrive at
(\ref{eq:jpdf_beta_alpha}).
\end{proof}

It is straightforward to obtain the joint PDF of $s$ and $r$ from
(\ref{eq:jpdf_beta_alpha}), according to these one-to-one mappings
$s=\sqrt{\beta}$ and $r=\sqrt{\alpha}$.

The joint PDF in (\ref{eq:jpdf_beta_alpha}) includes many existing
PDF's as special cases.
\begin{itemize}
\item By integration over $\beta$, (\ref{eq:jpdf_beta_alpha}) reduces to the
marginal PDF
\begin{equation}\label{eq:PDF_alpha}
p(\alpha)=\frac{1}{m}\sum_{k=0}^{m-1}\frac{k!}{(k+\nu)!}
\left[L_k^\nu(\alpha)\right]^2\alpha^\nu e^{-\alpha},
\end{equation} which is the same as the PDF presented
in \cite{IEEE_sw27:Telatar99_MIMO_Cap}. When $m=1$,
(\ref{eq:PDF_alpha}) further reduces to
\begin{equation}\label{eq:PDF_alpha_m=1}
p(\alpha)=\frac{1}{(n-1)!}\alpha^{n-1} e^{-\alpha},
\end{equation} which is the $\chi^2$ distribution
with $2n$ degrees of freedom\cite[(2.32)]{IEEE_sw27:SimonBook02}.

\item With $m=1$, (\ref{eq:jpdf_beta_alpha}) reduces
to\cite{IEEE_sw27:Wang06_Corr_Coeff_MRC},
\begin{equation}\label{eq:jpdf_alpha_beta_MRC}
p(\alpha,\beta)=\frac{(\alpha
\beta)^\frac{n-1}{2}\exp\left(\!-\frac{\alpha+\beta}{1-|\rho|^2}\!\right)
I_{n-1}\!\left(\frac{2|\rho|\sqrt{\alpha\beta}}{1-|\rho|^2}\right)}
{(n-1)!\left(1-|\rho|^2\right)|\rho|^{n-1}}.
\end{equation} Furthermore, when $n=1$,
(\ref{eq:jpdf_alpha_beta_MRC}) simplifies to
\begin{equation}\label{eq:jpdf_alpha_beta_SISO}
p(\alpha,
\beta)=\frac{1}{1-|\rho|^2}\exp\left(\!-\frac{\alpha+\beta}
{1-|\rho|^2}\!\right)
I_0\!\left(\frac{2|\rho|\sqrt{\alpha\beta}}{1-|\rho|^2}\right),
\end{equation} which is identical to
(8-103)\cite[pp. 163]{IEEE_sw27:DavenportBook87}, after two
one-to-one nonlinear transformations.
\end{itemize}

For the application discussed in section \ref{sec:Applications}, we
need the joint marginal PDF of $\phi$ and $\varphi$,
$p(\phi,\varphi)$, where $\phi$ and $\varphi$ are randomly selected
from $\{\alpha_k\}_{k=1}^m$, $m\geq2$. Using the technique in
\cite{IEEE_sw27:Mehta94_Two_Coupled_Hermitian_Matrices,
IEEE_sw27:Mehta98_Coupled_Hermitian_Matrix_Chain}, we have the
following theorem.
\begin{theorem}\label{th:jpdf_phi_varphi}
If $\phi$ and $\varphi$ are randomly selected from
$\{\alpha_k\}_{k=1}^m$, their joint PDF is given by
\begin{equation}\label{eq:jpdf_phi_varphi}
p(\phi,\varphi)\!=\!\frac{(\phi\varphi)^\nu
e^{-(\phi+\varphi)}}{m(m-1)}\!\sum_{\begin{subarray}{c}k,l=0\\k\neq
l\end{subarray}}^{m-1}\!\frac{k!l!}{(k+\nu)!(l+\nu)!}
\!\left\{\!\left[L_k^\nu(\phi)L_l^\nu(\varphi)\right]^2
\!-\!L_k^\nu(\phi)L_l^\nu(\phi)L_k^\nu(\varphi)L_l^\nu(\varphi)\!\right\}\!.
\end{equation}
\end{theorem}
\begin{proof}
According to (1.6) and (2.14) in
\cite{IEEE_sw27:Mehta94_Two_Coupled_Hermitian_Matrices} we have
\begin{equation}
p(\phi,\varphi)=\frac{1}{m(m-1)}\det\!\left|\!\!\!
                \begin{array}{cc}
                  K(\phi, \phi) & K(\phi, \varphi) \\
                  K(\varphi, \phi) & K(\varphi, \varphi) \\
                \end{array}
              \!\!\!\right|,
\end{equation} where
$K(x_1,x_2)=\sum_{k=0}^{m-1}P_k(x_1)\overline{Q}_k(x_2)$. With
$P_k(x_1)$ in (\ref{eq:Bi-orthogonal_Poly_P}) and
$\overline{Q}_k(x_2)$ in (\ref{eq:Weighted_Bi-orthogonal_Poly_Q}),
we obtain (\ref{eq:jpdf_phi_varphi}) after some simple algebraic
manipulations.
\end{proof}

\section{Application to Wireless MIMO Communication
Systems}\label{sec:Applications} For an $N_{\!R}\times N_{\!T}$ MIMO
time-varying Rayleigh flat fading channel\cite{IEEE_sw27:TseBook}
with $N_{\!T}$ transmitters and $N_{\!R}$ receivers, the channel
impulse response at time instant $t$ is given by
\begin{equation}\label{eq:H(t)}
\mathbf{H}(t)=\begin{bmatrix}
   h_{1,1}(t) &\cdots  & h_{1,N_{\!T}}(t) \\
   \vdots & \ddots & \vdots \\
   h_{N_{\!R},1}(t) & \cdots & h_{N_{\!R},N_{\!T}}(t) \\
 \end{bmatrix}\!.
\end{equation} We assume all the $N_{\!R}N_{\!T}$ sub-channels in the MIMO system,
$\left\{h_{i,j}(t)\right\}_{(i=1,j=1)}^{(N_{\!R},N_{\!T})}$ are
i.i.d., with the same temporal correlation coefficient, i.e.,
\begin{equation}\label{eq:Channel_Corr_Assumption}
\mathbb{E}[h_{ij}(t)h_{pq}^\star(t-\tau)]=\delta_{ip}\delta_{jq}
\rho_h(\tau),
\end{equation}
where $\rho_h(\tau)=J_0(2\pi
f_{\!D}\tau)$\cite{IEEE_sw27:JakesBook94} in isotropic scattering
environments\footnote{In the non-isotropic scattering environment,
$\rho_h(\tau)$, in general, is a complex-value
function\cite{IEEE_sw27:Wang05_LCR_AFD_SISO_Capacity,
IEEE_sw27:Wang05_IT}, and $|\rho_h(\tau)|$ indicates its amplitude
at the time delay $\tau$.}, with $J_0(x)=I_0(-\jmath x)$\cite[pp.
961, 8.406.3]{IEEE_sw27:RyzhikBook_5th} and $f_{\!D}$ is the maximum
Doppler frequency shift.

We set $n=\max(N_{\!R},N_{\!T})$ and $m=\min(N_{\!R},N_{\!T})$.
According to (\ref{eq:diagonalization_H}), $\mathbf{H}(t)$ can be
diagonalized as
\begin{equation}\label{eq:H(t)_SVD}
\mathbf{H}(t)=\mathbf{U}(t)\mathbf{S}(t)\mathbf{V}^\dag(t),
\end{equation} where
$\mathbf{S}(t)=\begin{bmatrix}\diag(\mathbf{s}(t))\
\mathbf{0}\end{bmatrix}$ with $\mathbf{s}(t)=[s_1(t), s_2(t),
\cdots, s_m(t)]$ for $N_{\!R}\leq N_{\!T}$, and
$\mathbf{S}(t)=\begin{bmatrix}\diag(\mathbf{s}(t))\\
\mathbf{0}\end{bmatrix}$ for $N_{\!R}>N_{\!T}$. Therefore the MIMO
channel, $\mathbf{H}(t)$, is decomposed to $m$ identically
distributed \emph{eigen}-channels $\lambda_k(t)=s_k^2(t)$,
$k\in[1,m]$, by SVD.

In wireless MIMO communication systems, we are interested in the
correlation coefficient between any two \emph{eigen}-channels, which
is defined by
\begin{equation}
\rho_{k,l}(\tau)=\frac{\mathbb{E}\left[\lambda_k(t)
\lambda_l(t-\tau)\right]-\mathbb{E}\left[\lambda_k(t)\right]
\mathbb{E}\left[\lambda_l(t)\right]}
{\sqrt{\mathbb{E}\left[\lambda_k^2(t)\right]
-\left\{\mathbb{E}\left[\lambda_k(t)\right]\right\}^2}
\sqrt{\mathbb{E}\left[\lambda_l^2(t)\right]
-\left\{\mathbb{E}\left[\lambda_l(t)\right]\right\}^2}}.
\end{equation}

For simplicity, in this paper we only consider a $2\times2$ MIMO
system, $N_{\!R}=N_{\!T}=2$, where the correlation coefficient,
$\rho_{k,l}(\tau)$, can be shown to be
\begin{equation}\label{eq:Coeff_Eigen-Channel}
\rho_{k,l}(\tau)=\begin{cases} 1-\frac{3}{2}\left(1-\delta_{kl}\right),&\tau=0,\\
\frac{|\rho_h(\tau)|^2}{4}=\frac{J_0^2\left(2\pi
f_{\!D}\tau\right)}{4}, &\tau\neq0,
\end{cases}, k,l=1,2,
\end{equation}
with $J^2_0(\cdot)=\left[J_0(\cdot)\right]^2$. To derive
(\ref{eq:Coeff_Eigen-Channel}), we note that for $\tau=0$ and $k=l$,
$\rho_{k,l}(0)=1$ because of the definition of the correlation
coefficient. Since $m=2$, for any \emph{eigen}-channel at the time
instant $t$, it is easy to show that the mean value of
$\lambda_k(t)$ is $\mathbb{E}\left[\lambda_k(t)\right]=2$, $\forall
k$, and the second moment of $\lambda_k(t)$ is
$\mathbb{E}\left[\lambda^2_k(t)\right]=8$, $\forall k$, using the
PDF in (\ref{eq:PDF_alpha}). For $\tau=0$ and $k\neq l$, we obtain
$\mathbb{E}\left[\lambda_k(t) \lambda_l(t)\right]=2$ by
(\ref{eq:jpdf_phi_varphi}), hence $\rho_{k,l}(0)=-\frac{1}{2}$,
$\forall k\neq l$. For $\tau\neq0$ and $\forall k,l$, it is not
difficult to get $\mathbb{E}\left[\lambda_k(t)
\lambda_l(t-\tau)\right]=4+|\rho_h(\tau)|^2$ using
(\ref{eq:jpdf_beta_alpha}), therefore we have the second line in
(\ref{eq:Coeff_Eigen-Channel}).

\begin{figure}[tp]
\centering
\includegraphics[width=.9\linewidth]
{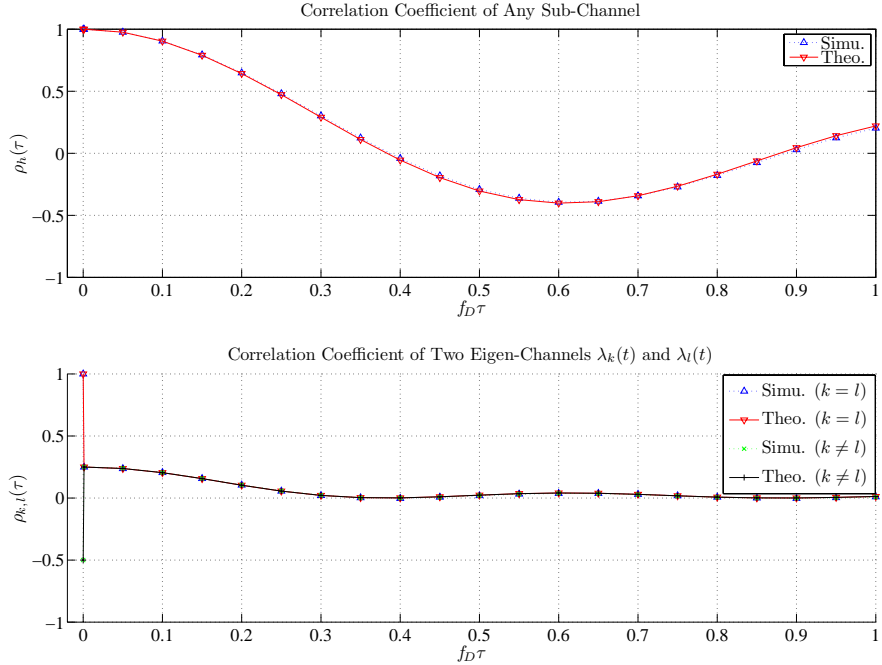} \caption{The channel correlation
coefficient, $\rho_h(\tau)$, and correlation coefficient of any two
eigen-channels, $\rho_{k,l}(\tau)$, in a $2\times2$ MIMO system with
Clarke's correlation model. Note that the sampling period,
$T_{\!s}$, is $\frac{1}{1000f_{\!D}}$ in Monte Carlo simulations,
therefore the first non-zero $\tau$ is $T_{\!s}$, i.e.,
$\frac{1}{1000f_{\!D}}$, which corresponds to
$f_{\!D}\tau=\frac{1}{1000}$ in the horizontal
axis.}\label{fig:CrossCoeff_Iso}
\end{figure}
Monte Carlo simulations are performed to verify the result in
(\ref{eq:Coeff_Eigen-Channel}). In all simulations\footnote{The
spectral method\cite{IEEE_sw27:Acolatse03} is used to generate the
MIMO channels.}, the maximum Doppler frequency $f_{\!D}$ is set to
$1$ Hz, and the sampling period, $T_{\!s}$, is equal to
$\frac{1}{1000f_{\!D}}$. The simulation results are shown in {\sc
Fig}. \ref{fig:CrossCoeff_Iso}, where the upper figure shows the
channel correlation coefficient $\rho_h(\tau)=J_0\left(2\pi
f_{\!D}\tau\right)$, Clarke's correlation model, whereas the lower
figure presents the correlation coefficient between any two
\emph{eigen}-channels or for any individual eigen-channel, Eq.
(\ref{eq:Coeff_Eigen-Channel}). Since $J_0(2\pi f_{\!D}\tau)$ is an
even function of $\tau$, the correlation coefficients  are plotted
for $\tau\geq0$. In all figures, ``Simu.'' indicates the curve is
obtained by Monte Carlo simulations, whereas ``Theo.'' means
theoretical. From {\sc Fig}. \ref{fig:CrossCoeff_Iso} we can
conclude that the new theoretical result in
(\ref{eq:Coeff_Eigen-Channel}) is confirmed by simulation very well.

\section{Conclusion}\label{sec:Conclusion}
In this paper, the joint distribution of singular values of two
correlated rectangular complex Gaussian random matrices is derived,
as well as the joint marginal distribution. The derived
distributions play an important role in the analysis and design of
wireless MIMO communication systems. As an example, the correlation
coefficient of any two \emph{eigen}-channels of a $2\times2$ MIMO
system is obtained and verified by the Monte Carlo simulations in
this paper.
\bibliographystyle{siam}
\bibliography{IEEEabrv,IEEE_sw27}
\end{document}